\documentclass{pja00}%
\usepackage{amsmath}
\usepackage{amsfonts}
\usepackage{amssymb}
\usepackage{graphicx}%
\newtheorem{Th}{Theorem}[section]
\newtheorem{Def}[Th]{Definition}
\newtheorem{Co}[Th]{Corollary}
\newtheorem{Lem}[Th]{Lemma}

\newtheorem{Pro}[Th]{Proposition}

\newcommand{\calS}{{\cal S}}

\providecommand{\U}[1]{\protect\rule{.1in}{.1in}}
\title{The Shi arrangement of the type $D_\ell$}
\Author{1}{Ruimei}{Gao}
\Author{1}{Donghe}{Pei}
\Author{2}{Hiroaki}{Terao}
\affiliation{1}{School of Mathematics and Statistics, Northeast Normal
University, Changchun 130024, P.R.China.
Partially supported by NSF (10871035) of China}
\affiliation{2}{Department of Mathematics, Hokkaido University, Sapporo, Hokkaido 060-0810, Japan. Partially supported by JSPS KAKENHI (21340001)}
\Subject[2010]{32S22;05E15}
\begin{document}
%
\maketitle
%

\begin{abstract}
In this paper, we give a basis for the derivation module of the cone over 
the Shi arrangement of the type $D_\ell$ explicitly.
\end{abstract}
\section{Introduction.}
Let $V$ be an $\ell$-dimensional vector space. An 
{\it
affine arrangement of hyperplanes}
 $\mathcal{A}$ is a finite collection of affine hyperplanes in $V$. If every hyperplane $H\in\mathcal{A}$ goes through the origin, then $\mathcal{A}$ is called to be {\it central}. When $\mathcal{A}$ is central, for each $H\in\mathcal{A}$, choose $\alpha_H\in V^*$ with $\ker(\alpha_H)=H$. Let $S$ be the algebra of polynomial functions on $V$ and let $\mathrm{Der}_{S}$ be the module of derivations 
\begin{align*}
&\mathrm{Der}_{S}:=\{\theta:S\rightarrow S\mid\theta(fg)=f\theta(g)+g\theta(f),f,g\in S,\\
&~~~~~~~~~~~~~~~~~~~~~~~~~~~~\theta~\mbox{is}~\mathbb{R}\mbox{-linear}\}.
\end{align*}
For a central arrangement $\mathcal{A}$, recall 
\begin{align*}
D(\mathcal{A}):=\{\theta\in\mathrm{Der}_{S}\mid\theta(\alpha_{H})
\in\alpha_{H}S
\mbox{~for all} ~H\in\mathcal{A}\}.
\end{align*}
We say that $\mathcal{A}$ is a {\it free arrangement} if $D(\mathcal{A})$ is a free $S$-module. The freeness was defined in \cite{Ter1}. The \linebreak Factorization Theorem\cite{Ter2} states that, for any free
arrangement $\mathcal{A}$, the characteristic polynomial of $\mathcal{A}$ factors completely over the  integers. 
%

Let $E=\mathbb{R}^\ell$ be an $\ell$-dimensional Euclidean \linebreak space  with a coodinate system $x_1,\ldots,x_\ell$, and $\Phi$ be a crystallographic irreducible root system. Fix a positive root system $\Phi^{+}\subset\Phi$. For each positive root $\alpha\in\Phi^{+}$ and $k\in\mathbb{Z}$, we define an affine hyperplane
\begin{align*}
    H_{\alpha,k}:=\{v\in V\mid (\alpha, v)=k\}.
\end{align*}
In \cite{Shi1}, J.-Y. Shi introduced the {\it Shi arrangement} 
\[
{\mathcal S}(A_{\ell})
:=
\{H_{\alpha,k}\mid\alpha\in\Phi^{+},\
0\leq k \leq 1\}
\]
when the root system is of the type $A_{\ell}$. 
This definition was later extended to the
{\it generalized Shi arrangement} (e.g., \cite{Edel})
\begin{align*}
\calS(\Phi):=
\{H_{\alpha,k}\mid\alpha\in\Phi^{+},\
0\leq k \leq 1\}.
\end{align*}
Embed $E$ into $V=\mathbb{R}^{\ell+1}$ by adding a new
coordinate $z$ such that $E$ is defined by the equation $z = 1$ in $V$.
Then, as in \cite{OT}, we have the cone
$
\mathbf{c}
\calS(\Phi)
$ 
of
$
\calS(\Phi)
$ 
\begin{align*}
\mathbf{c}
\calS(\Phi)
:=\{\mathbf{c}H_{\alpha,k} 
\mid\alpha\in\Phi^{+},\
0\leq k \leq 1\}\cup\{\{z=0\}\}.
\end{align*}
In \cite{Yo}, M. Yoshinaga 
proved that the cone 
$\mathbf{c}
\calS(\Phi)
$ 
is a free arrangement with exponenets
$(1, h, \dots, h)$
($h$ appears $\ell$ times),
where $h$ is the Coxeter number of $\Phi.$ 
(He actually verified the conjecture 
by P. Edelman and V. Reiner
in \cite{Edel}, which is far more general.)
He proved the freeness without finding a basis.

In \cite{Su1}, for the first time, the authors gave an 
explicit construction of
a basis for 
$D(\mathbf{c}\calS(A_{\ell}))$.
Then
D. Suyama constructed bases 
for 
$D(\mathbf{c}\calS(B_{\ell}))$
and
$D(\mathbf{c}\calS(C_{\ell}))$
in \cite{Su2}.
In this paper, we will give an explicit 
construction of a basis
for 
$D(\mathbf{c}\cal S(D_{\ell}))$.
A defining polynomial of
the cone over the
Shi arrangement of the type $D_{\ell}$  
is given by
\begin{align*}
Q := z \prod_{1\leq s<t\leq \ell} \prod_{\epsilon\in\{-1, 1\}} 
(x_s + \epsilon x_t - z)
(x_s + \epsilon x_t).
\end{align*}
Note that the number of hyperplanes in $\mathbf{c}\cal S(D_{\ell})$
is equal to $2\ell(\ell-1)+1$. 
Our construction is similar to the construction
in the case of the type $B_{\ell} $. 
The essential ingredients of the recipe are
 the Bernoulli polynomials and their \linebreak relatives.

\section{The basis construction.}

\bigskip

\begin{Pro}
\label{Prop2.1}
For $(p,q)\in\mathbb{Z}_{\geq -1}\times\mathbb{Z}_{\geq 0}$, \linebreak consider the following two conditions for a rational function $B_{p,q}(x)$:
\begin{align*}
&1.~B_{p,q}(x+1)-B_{p,q}(x)\\
&=\frac{(x+1)^p-(-x)^p}{(x+1)-(-x)}(x+1)^q(-x)^q,\\
&2.~B_{p,q}(-x)=-B_{p,q}(x).
\end{align*}
Then such a rational function $B_{p, q}(x) $ 
 uniquely exists.
Morever, the $B_{p, q} (x)$ is a polynomial unless $(p, q) = (-1, 0)$
and $B_{-1, 0}(x) = -(1/x) $. 
\end{Pro}

\medskip

\begin{proof}
Suppose $(p, q) \neq (-1, 0)$.  Since
the right hand side of the first condition is a polynomial in $x$,
there exists a polynomial $B_{p, q} (x)$ satisfying the first condition.
Note that $B_{p,q}(x)$ is unique up to a constant term.
Define a polynomial
$
F(x)=B_{p,q}(x)+B_{p,q}(-x)$.
Since
\begin{align*}
&~~~~B_{p,q}(-x)-B_{p,q}(-x-1)\\
&=\frac{(-x)^p-(x+1)^p}{(-x)-(x+1)}(-x)^q(x+1)^q\\
&=\frac{(x+1)^p-(-x)^p}{(x+1)-(-x)}(x+1)^q(-x)^q\\
&=B_{p,q}(x+1)-B_{p,q}(x),
\end{align*}
we have
$
F(x+1)=F(x)
$
for any $x$.  Therefore
$F(x) $ is a constant function.
Then the polynomial 
$B_{p, q}(x) -  \left(F(0)/2\right)$ is the unique solution
satisfying the both conditions.
Next
we suppose $(p,q)=(-1, 0)$.  Then we compute
\begin{align*}
&~~~~B_{-1, 0}(x+1)
-
B_{-1, 0}(x)\\
&=
\frac{(x+1)^{-1} -(-x)^{-1} }{(x+1)-(-x)}
=-\frac{1}{x+1}+\frac{1}{x}.
\end{align*}
Thus
$
B_{-1, 0}(x)=-(1/x)
$
is the unique solution
satisfying the both conditions.
\end{proof}

\medskip

\begin{Def} Define a rational function $\overline{B}_{p,q}(x,z)$ in $x$ and $z$ by
$$\overline{B}_{p,q}(x,z):=z^{p+2q}B_{p,q}(x/z).$$
\end{Def}
Then 
$\overline{B}_{p,q}(x,z)$ is a homogeneous polynomial of degree
$p+2q$ except the two cases:
$\overline{B}_{-1, 0} (x, z) = -(1/x)$
and
$\overline{B}_{0, q} (x, z) = 0$.

For a set $I := \{y_{1} , \dots, y_{m} \}$ of variables,
let $$
\sigma^{I}_{n}
:=
\sigma_{n} (y_{1} , \dots, y_{m} ),
\,\,
\tau^{I}_{2n}
:=
\sigma_{n} (y^{2}_{1} , \dots, y^{2}_{m} ),
$$ 
where $\sigma_{n} $ stands for
the elementary symmetric function of
degree $n$.

\medskip

\begin{Def} 
Define derivations
%
\begin{align*}
\varphi_j
&:=(x_j-x_{j+1}-z)
\sum\limits_{i=1}^\ell
\sum\limits_{\substack{K_{1}\cup K_{2} \subseteq J\\
K_{1} \cap K_{2} = \emptyset}}
\left(\prod K_{1} \right)
\left(\prod K_{2}\right)^{2} \\
&
(-z)^{|K_{1}|} 
\sum\limits_{\substack{0\leq n_1 \leq |J_{1}|\\ 0\leq n_2\leq |J_{2} |}}(-1)^{n_1+n_2}\sigma_{n_{1} }^{J_{1} } 
\tau_{2 n_{2}}^{J_{2} } 
\overline{B}_{k, k_{0}}(x_{i}, z)
\frac{\partial}{\partial x_{i}}
\end{align*}
for $j = 1,\dots, \ell-1$
and
\begin{align*}
\varphi_{\ell} 
&:=\sum\limits_{i=1}^\ell
\sum\limits_{\substack{K_{1}\cup K_{2} \subseteq J  \\
K_{1} \cap K_{2} = \emptyset}}
\left(\prod K_{1} \right)
\left(\prod K_{2}\right)^{2} 
(-z)^{|K_{1}|}\\
&
~~~~~~~~~~~
(-x_{\ell}) 
\overline{B}_{-1, k_{0}}(x_{i}, z)
\frac{\partial}{\partial x_{i}}
\end{align*}
for $j=\ell$,
where
\begin{align*} 
J &:= \{x_{1}, \dots, x_{j-1}\},\,
J_{1} := \{x_{j}, x_{j+1}\},\\
J_{2} &:= \{x_{j+2}, \dots, x_{\ell}\},\\
\prod K_{p} &:=\prod_{x_{i} \in K_{p}} x_{i}\,\,\,(p=1, 2),\\
k_{0} &:= |J \setminus (K_{1} \cup K_{2}) |\geq 0, \\
k 
&:= 
(|J_{1} |-n_{1}) +2(|J_{2} |-n_{2}) -1\geq -1.
\end{align*} 
\end{Def}

\medskip

Note
that
$\varphi_{j} (z)=0\,\,(1\leq j\leq\ell)$.
In the rest of the paper, we will give a proof of the following theorem:

\medskip

\begin{Th}
\label{main}  
The derivations 
$
\varphi_{1}, \dots, \varphi_{\ell},   
$
together with the Euler derivation 
$$\theta_{E} := 
z \frac{\partial}{\partial z} 
+ \sum_{i=1}^{\ell} x_{i} \frac{\partial}
{\partial x_{i}},
$$ 
form a basis for $D({\mathbf c}{\mathcal S}(D_{\ell}) )$.
\end{Th}

\medskip

Note 
that
$\theta_{E} (x_{i} )=x_{i} 
\,\,(1\leq i\leq\ell)$
and
$\theta_{E} (z)=z$.

\begin{Lem}
\label{Lemma2.5} 
Let
$1\leq i\leq \ell$
and
$1\leq j\leq \ell$.
Suppose $\varphi_{j}(x_{i})$ is nonzero.
Then $\varphi_{j}(x_{i} )$ is 
a homogeneous
polynomial
of degree
$2(\ell-1)$.
 \end{Lem}

\medskip
 
\begin{proof}
Define
\begin{align*}
F_{ij}&:=
(x_{j}-x_{j+1}-z)
\left(\prod K_{1} \right)
\left(\prod K_{2}\right)^{2}
z^{|K_{1}|} \\
&~~~~~~
\sigma_{n_{1} }^{J_{1} }
\tau_{2 n_{2}}^{J_{2} } 
\overline{B}_{k, k_{0}}(x_{i}, z)~~~~~~~~~(1\leq j\leq \ell-1),\\
F_{i\ell}
&:=
\left(\prod K_{1} \right)
\left(\prod K_{2}\right)^{2}
z^{|K_{1}|} 
x_{\ell} 
\overline{B}_{-1, k_{0}}(x_{i}, z)
\end{align*}
when $K_{1} , K_{2} , n_{1} , n_{2} $ are fixed.
Then $\varphi_{j} (x_{i} )$ is a linear combination of
the $F_{ij} $'s over ${\mathbb R}$.  

Note that
$\overline{B}_{k, k_{0}} (x_{i}, z)$
is a polynomial unless 
$(k, k_{0}) = (-1, 0)$.

Assume that $1\leq j\leq \ell-1$ and $(k, k_{0} )=(-1, 0)$.
Then
$J = K_{1} \cup K_{2} $,
$n_{1}  = |J_{1}|$,
$n_{2}  = |J_{2}|$,
and 
$\overline{B}_{-1, 0} (x_{i}, z) =-1/x_{i}$.
Therefore
each 
$F_{i j}$ 
is a polynomial.
Thus 
$\varphi_{j}(x_{i})$
is a nonzero polynomial and
there exists a nonzero polynomial $F_{ij} $. 
Compute
\begin{align*}
&~~~\deg \varphi_{j}(x_{i} )=\deg F_{ij} \\
&=
1+ |K_{1}| + 2|K_{2}| 
+|K_{1}|
+n_{1} +2n_{2} 
\\
&
~~~~~~~~~~+
\deg \overline{B}_{k, k_{0}}(x_{i}, z)\\
&=
1+ 2|K_{1}| + 2|K_{2}| 
+
n_{1} +2n_{2} +
(2 k_{0}+k)\\
&=
1+ 2|K_{1}| + 2|K_{2}| 
+n_{1} +2n_{2}\\ 
&~~~+
2 (|J|-|K_{1}| -|K_{2}|)
+|J_{1} | -n_{1}\\
&~~~+2(|J_{2}|- n_{2}) -1\\
&=
2(|J|+|J_{1} |+|J_{2} |)-|J_{1} |=
2 \ell-2.
\end{align*}

Next consider
$\varphi_{\ell}(x_{i})$.
If
$k_{0} = 0$, then
$J = K_{1} \cup K_{2}$.
Therefore
each 
$F_{i \ell}$
is a polynomial.
Thus so is
$\varphi_{\ell}(x_{i})$. 
Compute
\begin{align*}
&~~~\deg \varphi_{\ell}(x_{i} )\\
&=
|K_{1}| + 2|K_{2}| 
+|K_{1}|+
1+
\deg \overline{B}_{-1, k_{0}}(x_{i}, z)\\
&=
2(|K_{1}| + |K_{2}|) 
+1+
(2 k_{0}-1)\\
&=
2(|K_{1}|+|K_{2} |+k_{0})=
2(\ell-1).
\end{align*}
\end{proof}

Let $<$ denote the {\it
pure lexicographic order} of monomials
with respect to the total order
\[
x_{1} > x_{2} >\dots > x_{\ell} > z.
\]
When $f\in S={\mathbb C}[x_{1} , x_{2} , \dots, x_{\ell}, z ] $
 is a nonzero polynomial, let 
$\mathrm{in} (f)$ denote the {\it
initial monomial} (e.g., see \cite{Her}) of $f$
with respect to the order
$<$.
%

\medskip

\begin{Pro}
\label{in}
Suppose $\varphi_{j} (x_{i} )$ is nonzero.
Then

(1)
$\mathrm{in}(\varphi_j(x_i))
\leq
x_1^2\cdots x_{i-1}^2x_i^{2\ell-2i}$,

(2)
$
\mathrm{in}(\varphi_j(x_i))<
x_1^2\cdots x_{i-1}^2x_i^{2\ell-2i}$
for $i < j,$ 

(3)
$
\mathrm{in}(\varphi_i(x_i))=x_1^2\cdots x_{i-1}^2x_i^{2\ell-2i}$
for
$1\leq i\leq \ell.
$
\end{Pro}

\smallskip

\begin{proof}
Recall 
$F_{ij}
\,(1\leq j\leq \ell-1)
$
and
$F_{i\ell}
$
from the proof of Lemma
\ref{Lemma2.5} 
when
$K_1,K_2,n_1,n_2$ are fixed.
Let $\deg^{(x_i)}f$ denote the degree of $f$
with respect to $x_{i} $ when $f\neq 0$.

(1)
Since,
 for every
nonzero $F_{ij}$,
we obtain
\begin{align*}
\deg^{(x_{p})}F_{ij}\leq 2
\,\,(1\leq p < i),
\,\,\,\,\,\,
\deg(F_{ij})=2\ell-2.
\end{align*}
Hence
we may conclude 
\begin{align*}
\mathrm{in}(F_{ij})\leq 
x_1^2\cdots x_{i-1}^2x_i^{2\ell-2i}
\end{align*}
and thus\begin{align*}
&\mathrm{in}(\varphi_j(x_i))
\leq\max\{\mathrm{in}(F_{ij})\}
\leq x_1^2\cdots x_{i-1}^2x_i^{2\ell-2i}.
\end{align*}

(2)
Suppose $i<j<\ell$. 
Since $x_{i} > x_{j} > z$, one has 
\begin{align*}
&~~~\mathrm{in} (\sigma_{n_{1} }^{J_{1}}
\tau_{2 n_{2}}^{J_{2}}\overline{B}_{k,k_0}(x_i,z))\\
&\leq x_{i}^{n_{1} +
2 n_{2} +2 k_{0}+k }
=
x_{i}^{2\ell-2j+2k_0-1}
\end{align*}
when $\overline{B}_{k,k_0}(x_i,z)$ is nonzero.
The equality holds if and only if $n_1=n_2=0$.

Suppose that $F_{ij} $ is nonzero.
For $1\leq i<j\leq \ell-1$, we have
\begin{align*}
&\mathrm{in}(F_{ij})\\
=&\mathrm{in}(x_j-x_{j+1}-z)\mathrm{in}\big((\prod K_1)(\prod K_2)^2(-z)^{|K_1|}\big)\\
&\mathrm{in}\big(\sigma_{n_{1} }^{J_{1}}
\tau_{2 n_{2}}^{J_{2}}
\overline{B}_{k,k_{0}}(x_{i}, z)\big)\\
\leq& x_j\,
\mathrm{in}\big((\prod K_1)(\prod K_2)^2(-z)^{|K_1|}\big)x_i^{2\ell-2j+2k_0-1}\\
=& x_j\,\mathrm{in}\big((\prod K_1)(\prod K_2)^2(-z)^{|K_1|}x_i^{2k_0}\big)x_i^{2\ell-2j-1}\\
\leq&x_j(x_1^2\cdots x_{i-1}^2x_i^{2j-2i})x_i^{2\ell-2j-1}~~~(*)\\
=&x_1^2\cdots x_{i-1}^2x_i^{2\ell-2i-1}x_j
<x_1^2\cdots x_{i-1}^2x_i^{2\ell-2i}.
\end{align*}
Thus
\begin{align*}
\mathrm{in}(\varphi_j(x_i))<x_1^2\cdots x_{i-1}^2x_i^{2\ell-2i}.
\end{align*}
For $1\leq i<j=\ell$,
\begin{align*}
&\mathrm{in}(F_{i\ell})\\
=&x_\ell\,\mathrm{in}\big((\prod K_1)(\prod K_2)^2(-z)^{|K_1|}\big)\mathrm{in}\big(\overline{B}_{-1,k_{0}}(x_{i}, z)\big)\\
=& x_\ell\,\mathrm{in}\big((\prod K_1)(\prod K_2)^2(-z)^{|K_1|}\big)x_i^{2k_0-1}\\
=&x_\ell\,\mathrm{in}\big((\prod K_1)(\prod K_2)^2(-z)^{|K_1|}x_i^{2k_0}\big)x_i^{-1}\\
\leq&x_\ell(x_1^2\cdots x_{i-1}^2x_i^{2\ell-2i})x_i^{-1}~~~~(**)\\
=&x_1^2\cdots x_{i-1}^2x_i^{2\ell-2i-1}x_\ell
<x_1^2\cdots x_{i-1}^2x_i^{2\ell-2i}.
\end{align*}
This proves (2).

Now we only need to prove (3).
Let $i=j<\ell$ 
in $(*)$. Then the equality
\begin{align*}
&\mathrm{in}(F_{ii})
=x_1^2\cdots x_{i-1}^2x_i^{2\ell-2i}
\end{align*}
holds  if and only if 
\begin{align*}
K_1=\emptyset, \,K_2=J, \,n_1=n_2=k_{0} =0,
\, k=2\ell-2i-1
\end{align*}
because
the leading term of $\overline{B}_{2\ell-2i-1,0}(x_i,z)$
is equal to
$$\frac{x_{i}^{2\ell-2i-1}}{2\ell-2i-1}.$$
Next let $i=\ell$ 
in $(**)$.
Then the equality 
\begin{align*}
&\mathrm{in}(F_{\ell\ell})
=x_1^2\cdots x_{\ell-1}^2
\end{align*}
 holds
if and only if $$
K_1=\emptyset, \,K_2=J=\{x_{1} ,\dots,x_{\ell-1}\}, \, k_{0} = 0.$$ 
Therefore,
for $1\leq i\leq \ell$,
\begin{align*}
\mathrm{in}(\varphi_i(x_i))=x_1^2\cdots x_{i-1}^2x_i^{2\ell-2i}.
\end{align*}
\end{proof}

From Proposition \ref{in}, 
we immediately obtain the following Corollary:

\medskip

\begin{Co}
\label{Coin} 
(1)
$$
\mathrm{in}(\det\big[\varphi_{j}(x_i)\big])
=\prod_{i=1}^\ell\mathrm{in}(\varphi_i(x_i))
=
\prod_{i=1}^{\ell-1} x_i^{4(\ell-i)}.
$$

(2)
Moreover,
the leading term of
$\det\big[\varphi_{j}(x_i)\big]$
is equal to
\[
\frac{1}{(2\ell-3)!!}  
\prod_{i=1}^{\ell-1} x_i^{4(\ell-i)}.
\]

(3) In particular, 
$\det\big[\varphi_{j}(x_i)\big]$
does not vanish.
%
\end{Co}

\medskip

\medskip
Next, we will prove $\varphi_j\in D({\mathbb c}\mathcal{S}(D_\ell))$
for $1\leq j\leq \ell$.
We denote ${\mathbb c}\mathcal{S}(D_\ell)$ simply
by
$\mathcal S_{\ell}$ from now on.
Before the proof, we need the following two lemmas:

\medskip
\begin{Lem}
\label{Lem2.6} 
Fix $1\leq j\leq \ell-1$ and $\epsilon\in\{-1, 1\}$.  Then

\noindent 
(1)
\begin{multline*}
\prod_{x_{i}\in J }(x_i-x_s)(x_i-\epsilon x_t)
=\sum\limits_{\substack{K_{1}\cup K_{2} \subseteq J\\
K_{1} \cap K_{2} = \emptyset}}
\left(\prod
K_{1}
\right)\\
\times\left(\prod
K_{2}
\right)^{2}
%
[-(x_s+\epsilon x_t)]^{|K_1|}(\epsilon x_sx_t)^{k_0}.
\end{multline*}
(2)
\begin{multline*}
\sum\limits_{\substack{0\leq n_1 \leq |J_{1}|\\0\leq n_2\leq |J_{2}|}}
(-1)^{|J_{1}|+|J_{2}|-n_{1} -n_2}
\sigma_{n_1}^{J_{1}} \tau_{2n_2}^{J_{2}} 
(\epsilon  x_s)^{k+1}\\
=
\prod_{x_{i}\in J_{1}} (x_i- \epsilon x_s)
\prod_{x_{i}\in J_{2}}(x_i^2-x_s^2).
\end{multline*}
\end{Lem}

\vspace{3truemm}
\begin{proof}
(1) is easy because the left handside is equal to
\[
\prod\limits_{x_{i}\in J}(x_i^{2} -
(x_s+\epsilon x_t)x_{i}  + \epsilon x_{s} x_t).
\]
(2) The left handside is equal to
\begin{align*}
&
\sum\limits_{0\leq n_1 \leq |J_{1}|}
(-\epsilon x_s)^{|J_{1}|-n_1}
\sigma_{n_1}^{J_{1} } 
\sum\limits_{0\leq n_2\leq |J_{2} |}
(-x_s^2)^{|J_{2} |-n_2}
\tau_{2n_2}^{J_{2} } 
\end{align*}
which is equal to the right handside.
\end{proof}

\medskip
\begin{Lem}

\noindent
(1)The polynomial
$$
x_{s} \overline{B}_{k,k_{0}}(x_s,z)-x_{t} \overline{B}_{k,k_{0} }(x_t,z)
$$
is divisible by $x_{s}^{2} - x_{t}^{2},   $ 

\noindent
(2)For $\epsilon\in\{-1, 1\}$,
the polynomial
\begin{multline*}
(x_{s}-\epsilon x_{t})
\epsilon x_{s} x_{t} 
\left[
\overline{B}_{k,k_{0}}(x_s,z)+
\epsilon \overline{B}_{k,k_{0}}(x_t,z)\right]\\
- (x_s+\epsilon x_t)(\epsilon x_{s}x_t)^{k_0}
\left[\epsilon x_{t} x_s^{k+1}- x_{s}(\epsilon  x_t)^{k+1}\right]
\end{multline*}
is divisible by
$x_s+ \epsilon x_t-z$.
\label{Lem2.7} 
\end{Lem}
\vspace{3truemm}
\begin{proof}

(1) follows from the fact that $
-\overline{B}_{k,k_0}(x,z)
=\overline{B}_{k,k_0}(-x,z)
$
in Proposition \ref{Prop2.1}.

(2) follows from 
the following congruence relation of 
polynomials modulo
the ideal $(x_s+ 
\epsilon x_t-z)$: 
\begin{align*}&
(x_{s}-\epsilon x_{t})
\epsilon x_{s} x_{t} 
\left[
\overline{B}_{k,k_{0}}(x_s,z)+\epsilon 
\overline{B}_{k,k_{0}}(x_t,z)\right]\\
&=
(x_{s}-\epsilon x_{t})
\epsilon x_{s} x_{t} 
z^{k+2k_0}\big[{B}_{k,k_0}(\displaystyle\frac{x_s}{z})- 
{B}_{k,k_0}(\displaystyle\frac{ -\epsilon x_t}{z})\big]\\
&\equiv
(x_{s}-\epsilon x_{t})
\epsilon x_{s} x_{t} 
(x_s+ \epsilon x_t)^{k+2k_0}\\
&~~~~\big[{B}_{k,k_0}(\displaystyle\frac{x_s}{x_s+ \epsilon x_t})
- {B}_{k,k_0}
(\displaystyle\frac{ -\epsilon x_t}{x_s+ \epsilon x_t})
\big]\\
&=(x_{s}-\epsilon x_{t})
\epsilon x_{s} x_{t} 
(x_s+ \epsilon x_t)^{k+2k_0}\\
&~~~~\displaystyle\frac{\big(\displaystyle\frac{x_s}{x_s+ 
\epsilon x_t}\big)^k
-
\big(\displaystyle\frac{\epsilon x_t}{x_s+ \epsilon x_t}\big)^k}
{\big({\displaystyle\frac{x_s}{x_s+ \epsilon x_t}}\big)-
\big({\displaystyle\frac{ \epsilon x_t}{x_s+\epsilon x_t}}\big)}
({\displaystyle\frac{ \epsilon x_t}{x_s+ \epsilon x_t}})^{k_0}
({\displaystyle\frac{x_s}{x_s+ \epsilon x_t}})^{k_0}\\
&=(x_s+ \epsilon x_t)
(\epsilon  x_{s}x_t)^{k_0}
\left[\epsilon x_{t} x_s^{k+1}- x_{s}(\epsilon  x_t)^{k+1}\right].
\end{align*}
\end{proof}
\begin{Pro}Every
$\varphi_j$
lies in
$D(\mathcal{S}_\ell)$.
\label{Proposition2.10} 
\end{Pro}
\vspace{3truemm}
\begin{proof}
For
$1\leq j\leq \ell-1,1\leq s< t\leq \ell,$
and $\epsilon\in\{-1, 1\}$, 
by Lemma \ref{Lem2.7} and Lemma \ref{Lem2.6},
we have
the following congruence relation of 
polynomials modulo the ideal
$(x_{s} +\epsilon x_{t} -z)$: 
\begin{align*}
&
(x_{s}-\epsilon x_{t})
\epsilon x_{s} x_{t} 
\left[\varphi_j(x_s+\epsilon x_t-z)\right]\\
&=(x_j-x_{j+1}-z)
\sum\limits_{\substack{K_{1}\cup K_{2} \subseteq J\\
K_{1} \cap K_{2} = \emptyset}}
\left(\prod K_{1} \right)
\left(\prod K_{2}\right)^{2}\\
&~~~~~~~\times
(-z)^{|K_{1}|} 
\sum\limits_{\substack{0\leq n_1 \leq |J_{1}|\\ 0\leq n_2\leq |J_{2} |}}(-1)^{n_1+n_2}\sigma_{n_{1} }^{J_{1} } 
\tau_{2 n_{2}}^{J_{2} }\\
&~~~~~~~\times(x_{s}-\epsilon x_{t})
\epsilon x_{s} x_{t} 
[\overline{B}_{k,k_0}(x_s,z)+\epsilon \overline{B}_{k,k_0}(x_t,z)]
\\
&\equiv(x_j-x_{j+1}-z)
\left(x_{s}+\epsilon x_{t}\right)\\
&~\times
\sum\limits_{K_{1}, K_{2}}
\left(\prod K_{1} \right)
\left(\prod K_{2}\right)^{2} 
[-(x_s+\epsilon x_t)]^{|K_{1}|}(\epsilon x_{s}x_t)^{k_0}\\ 
&~\times
\sum\limits_{n_{1}, n_2}
(-1)^{n_{1} + n_{2} } 
\sigma_{n_1}^{J_{1} } \tau_{2n_2}^{J_{2} } 
\left[\epsilon x_{t} x_s^{k+1}- x_{s}(\epsilon  x_t)^{k+1}\right]\\
&=(x_j-x_{j+1}-z)\left(
x_s+\epsilon x_t\right)
\prod_{x_{i} \in J}(x_i-x_s)(x_i-\epsilon x_t)\\
&~~\times(-1)^{|J_{2}|} \bigg[
\epsilon x_{t}\prod_{x_{i}\in J_{1}  }  (x_{i}-x_s)
\prod_{x_{i}\in J_{2}}(x_i^2-x_s^2)
\\
&~~~~~~~~~~~~~~
-x_{s}
\prod_{x_{i}\in J_{1}  }  (x_{i}- \epsilon x_{t} )
\prod_{x_{i}\in J_{2}  }(x_i^2-x_t^2)
\bigg]\,\,\,\,\,\, (\dagger).
\end{align*}

\noindent
{\it Case 1.}
When $x_{s} \in J$, $(\dagger)=0.$  

\noindent
{\it Case 2.}
When 
$x_{s} \in J_{2} $
 and
$x_{t} \in J_{2} $, $(\dagger)=0.$  

\noindent
{\it Case 3.}
When 
$x_{s} \in J_{1} $
and
$x_{t} \in J_{2} $, $(\dagger)=0.$

\noindent
{\it Case 4.}
When
$x_{s} \in J_{1} $,
$x_{t} \in J_{1} $
and $\epsilon=1$, 
$(\dagger)=0$.  

\noindent
{\it Case 5.}
If
$x_{s} \in J_{1} $,
$x_{t} \in J_{1} $
and $\epsilon=-1$, then
$s=j<t=j+1$.  So $(\dagger)$ is divisible by $x_{s} +\epsilon x_{t} -z$. 

We also have
the following congruence relation of 
polynomials modulo the ideal
$(x_{s} +\epsilon x_{t} -z)$: 
\begin{align*}
&
(x_{s}-\epsilon x_{t})
\epsilon x_{s} x_{t} 
\left[\varphi_\ell(x_s+\epsilon x_t-z)\right]\\
&=
\sum\limits_{\substack{K_{1}\cup K_{2} \subseteq J\\
K_{1} \cap K_{2} = \emptyset}}
\left(\prod K_{1} \right)
\left(\prod K_{2}\right)^{2}
(-z)^{|K_{1}|} 
(-x_{\ell} )
\\
&
(x_{s}-\epsilon x_{t})
\epsilon x_{s} x_{t} 
[\overline{B}_{-1,k_0}(x_s,z)+\epsilon \overline{B}_{-1,k_0}(x_t,z)]
\\
&\equiv
\left(x_{s}+\epsilon x_{t}\right)
(-x_{\ell})
\left(\epsilon x_{t} - x_{s}\right)
\\
&
\sum\limits_{K_{1}, K_{2}}
\left(\prod K_{1} \right)
\left(\prod K_{2}\right)^{2} 
[-(x_s+\epsilon x_t)]^{|K_{1}|}(\epsilon x_{s}x_t)^{k_0}\\ 
&=
\left(
x_s^{2} - x_t^{2} \right)
x_{\ell} \prod_{x_{i} \in J}(x_i-x_s)(x_i-\epsilon x_t)
\,\,\,\,\,\,\,\,\,\,\,\, (\dagger\dagger).
\end{align*}
Since
$s<t\leq \ell$, we have $x_{s}\in J = \{x_{1}, \dots, x_{\ell-1}  \}$.
Thus $(\dagger\dagger)=0.$ 
\noindent
Therefore
$\varphi_j(x_s+\epsilon x_t-z)$ is divisible by 
$x_{s}+\epsilon x_{t}-z$ for 
$1\leq j\leq \ell, 1\leq s< t\leq \ell.$

For $1\leq j\leq \ell$, 
\begin{align*}
\varphi_j(x_s^{2} - x_t^{2} )
=
2 x_{s} \varphi_j(x_s) - 2 x_t 
\varphi_j(x_t)
\end{align*}
is divisible either by
$x_{s} \overline{B}_{k,k_{0}}(x_s,z)-x_{t} \overline{B}_{k,k_{0} }(x_t,z)
$
or
by
$x_{s} \overline{B}_{-1,k_{0}}(x_s,z)-x_{t} \overline{B}_{-1,k_{0} }(x_t,z)$,
we have
\[
\varphi_j(
x_s^{2} - x_t^{2} 
)
\equiv
0
\,\,
\mod
(x_s^{2} - x_t^{2})
\]
by Lemma \ref{Lem2.7} (1).
This implies $\varphi_{j} \in D(\cal S_{\ell} ).$
\end{proof}
\vspace{3truemm}
Applying Saito's lemma \cite{Sai} 
\cite[Theorem 4.19]{OT},
we complete our proof of Theorem \ref{main}
thanks to Lemma \ref{Lemma2.5}, Corollay \ref{Coin}  (3)
and
Proposition \ref{Proposition2.10}.  
Theorem \ref{main} implies that $\det[\varphi_{j} (x_{i})]$ 
is a nonzero multiple of 
$(Q/z)$.  By Corollary  \ref{Coin} (2)
one obtains

\begin{Co}
\begin{multline*} 
\det[\varphi_{j}(x_{i} ) ]\\
=
\frac{1}{(2\ell-3)!!}
\prod_{1\leq s<t\leq \ell} \prod_{\epsilon\in\{-1, 1\}} 
(x_s + \epsilon x_t - z)
(x_s + \epsilon x_t).
\end{multline*}
\end{Co}


\begin{thebibliography}{17}
\bibitem{AT}
T. Abe, H. Terao, The freeness of Shi-Catalan arrangements. {\em European J. Combin.}(to appear). arXiv:1012.5884v1.
\bibitem{Ath1}
Ch. Athanasiadis, Characteristic polynomials of subspace arrangements and finite fields. {\em Adv. in Math.}, {\bf 122} (1996), 193-233.
\bibitem{Ath2}
Ch. Athanasiadis, On free deformations of the braid arrangement. {\em European J. Combin.}, {\bf 19} (1998), 7-18.
\bibitem{Edel}
P. H. Edelman and V. Reiner, Free arrangements and rhombic tilings. {\em Discrete Comp. Geom.}, {\bf 15} (1996), 307-340.
\bibitem{Head} P. Headley, On a family of hyperplane arrangements related to the affine Weyl groups. {\em J. Algebraic Combin.}, {\bf 6} (1997), 331-338.
\bibitem{Her} J. Herzog and T. Hibi, {\em Monomial ideals.} Graduate Texts in Mathematics, Springer-Verlag, London, 2011.
\bibitem{OT}
P. Orlik and H. Terao, {\em Arrangements of hyperplanes.} Grundlehren der Mathematischen Wissenschaften, {\bf 300}, Springer-Verlag, Berlin, 1992.
\bibitem{PS}
A. Postnikov, R. P. Stanley, Deformations of Coxeter hyperplane arrangements. {\em J. Comb. Theory, Ser. A}, {\bf 91} (2000), 544-597.
\bibitem{Sai}
K. Saito, Theory of logarithmic differential forms and logarithmic vector fields. {\em J. Fac. Sci. Univ. Tokyo Sect. IA Math.}, {\bf 27} (1980), 265-291.
\bibitem{Shi1}
J.-Y. Shi, The Kazhdan-Lusztig cells in certain affine Weyl groups. {\em Lecture Notes in Math.}, {\bf 1179}, Springer-Verlag, 1986.
\bibitem{Shi2}
J.-Y. Shi, Sign types corresponding to an affine Weyl group. {\em J. Lond. Math. Soc.}, {\bf 35} (1987), 56-74.
\bibitem{ST}
L. Solomon and H. Terao, The double Coxeter arrangements. {\em
Comment. Math. Helv.}, {\bf 73} (1998), 237-258.
\bibitem{Su1}
D. Suyama, H. Terao, The Shi arrangements and the Bernoulli polynomials. arXiv:1103.3214v3.
\bibitem{Su2}
D. Suyama, On the Shi arrangements of types 
$B_{\ell} $,
$C_{\ell} $,
$F_{4} $
and $G_{2}. $ 
(in preparation).
 

\bibitem{Ter1}
H. Terao, Arrangements of hyperplanes and their freeness I, II. {\em J. Fac.
Sci. Univ. Tokyo}, {\bf 27} (1980), 293-320.
\bibitem{Ter2}
H. Terao, Generalized exponents of a free arrangement of hyperplanes
and Shephard-Todd-Brieskorn formula. {\em Invent. Math.}, {\bf 63} (1981), 159-179.
\bibitem{Ter3}
H. Terao, Multiderivations of Coxeter arrangements. {\em Invent. Math.}, {\bf 148} (2002), 659-674.
\bibitem{Yo}
M. Yoshinaga, Characterization of a free arrangement and conjecture of
Edelman and Reiner. {\em Invent. Math.}, {\bf 157} (2004), 449-454.
\end{thebibliography}
\end{document}